\documentclass[titlepage,12pt]{article}
\usepackage{amssymb}
\usepackage{latexsym}
\usepackage{pspic}

\newtheorem{theorem}{Theorem}[section]
\newtheorem{proposition}[theorem]{Proposition}
\newtheorem{corollary}[theorem]{Corollary}
\newtheorem{lemma}[theorem]{Lemma}

\def\su{\sum_{ \kappa } x_ {\kappa( v_{1} )}
x_{\kappa(v_{2})} \cdots x_{ \kappa ( v_{d})}}
\def\ind{\hspace{-.05 in}\uparrow}

\def\non{ Y_{G}}
\def\nond{  Y_{G - e}}
\def\nonc { Y_{G/e}}

\def\n{\noindent}

\def\comp{\circ}
\def\chr{{\cal X}_{G}(n)}

\def\v{v_{0}}
\newcommand{\ben}{\begin{enumerate}}
\newcommand{\een}{\end{enumerate}}
\newcommand{\ble}{\begin{lemma}}
\newcommand{\ele}{\end{lemma}}
\newcommand{\bth}{\begin{theorem}}
\renewcommand{\eth}{\end{theorem}}
\newcommand{\bpr}{\begin{proposition}}
\newcommand{\epr}{\end{proposition}}
\newcommand{\bco}{\begin{corollary}}
\newcommand{\eco}{\end{corollary}}
\newcommand{\bcon}{\begin{conj}}
\newcommand{\econ}{\end{conj}}
\newcommand{\bde}{\begin{defn}}
\newcommand{\ede}{\end{defn}}
\newcommand{\bex}{\begin{exa}}
\newcommand{\eex}{\end{exa}}
\newcommand{\barr}{\begin{array}}
\newcommand{\earr}{\end{array}}
\newcommand{\btab}{\begin{tabular}}
\newcommand{\etab}{\end{tabular}}
\newcommand{\beq}{\begin{equation}}
\newcommand{\eeq}{\end{equation}}
\newcommand{\bea}{\begin{eqnarray*}}
\newcommand{\eea}{\end{eqnarray*}}
\newcommand{\bce}{\begin{center}}
\newcommand{\ece}{\end{center}}
\newcommand{\bpi}{\begin{picture}}
\newcommand{\epi}{\end{picture}}
\newcommand{\bfi}{\begin{figure} \begin{center}}
\newcommand{\efi}{\end{center} \end{figure}}
\newcommand{\capt}{\caption}
\newcommand{\bsl}{\begin{slide}{}}
\newcommand{\esl}{\end{slide}}

\newcommand{\Qed}{\rule{1ex}{1ex} \medskip}

\newcommand{\emp}{\emptyset}

\newcommand{\setm}{\setminus}

\newcommand{\mh}{\hat{0}}

\newcommand{\mt}{\wedge}

\newcommand{\ree}[1]{(\ref{#1})}

\newcommand{\ra}{\rightarrow}

\newcommand{\ka}{\kappa}

\newcommand{\bx}{{\bf x}}

\newcommand{\bbC}{{\mathbb C}}

\newcommand{\bbP}{{\mathbb P}}

\newcommand{\cA}{{\cal A}}

\newcommand{\cD}{{\cal D}}

\newcommand{\ds}{\displaystyle}
\newcommand{\defn}{\stackrel{\rm def}{=}}

\newtheorem{definition}[theorem]{Definition}

\def\su{\sum_{ \kappa } x_ {\kappa( v_{1} )}
x_{\kappa(v_{2})} \cdots x_{ \kappa ( v_{d})}}
\def\ind{\hspace{-.05 in}\uparrow}

\def\non{ Y_{G}}
\def\nond{  Y_{G \setminus e}}
\def\nonc { Y_{G/e}}

\def\n{\noindent}

\def\comp{\circ}
\def\chr{{\cal X}_{G}(n)}

\begin{document}
\pagestyle{empty}
\title{Sinks in Acyclic Orientations of Graphs
}
\author{David D. Gebhard \\ 
Department of Mathematics\\
Campus Box 7281\\
Lenoir-Rhyne College\\
Hickory, NC 28603\\
gebhardd@lrc.edu\\[20pt]
Bruce E. Sagan\\
Department of Mathematics\\
Michigan State University\\
East Lansing, MI 48824-1027\\
sagan@math.msu.edu}

\date{\today \\[1in]
\begin{flushleft}
\end{flushleft}
       }
\maketitle

\begin{flushleft} Proposed running head:
\end{flushleft}
\begin{center}
Sinks in digraphs
\end{center}

Send proofs to:
\begin{center}
David D. Gebhard \\ 
Department of Mathematics\\
101C Minges Hall\\
Lenoir-Rhyne College\\
Hickory, NC 28601\\
gebhard@lrc.edu
\end{center}

\begin{abstract}
Greene and Zaslavsky proved that the number of acyclic orientations of
a graph with a unique sink at a given vertex is, up to sign, the
linear coefficient of the chromatic polynomial.  We give three new
proofs of this result using pure induction, noncommutative symmetric
functions, and an algorithmic bijection.
\end{abstract}
\pagestyle{plain}

\section{Acyclic Orientations}

Our primary focus will be a theorem of Greene and
Zaslavsky~\cite{gandz} concerning acyclic orientations of a graph and
its chromatic polynomial.  To state it we need some definitions.
For any undefined terms, we follow the terminology of Harary's
book~\cite{har:gt}.

Let $G$ be a finite graph with vertices $V=V(G)$ and edges $E=E(G)$.
We permit $E$ to contain loops and multiple edges.
An {\it orientation} of $G$ is a digraph formed by replacing $e\in E$
by one of the two possible directed arcs.  The orientation is
{\it acyclic} if it has no directed cycles.  We let $\cA(G)$ be the
set of acyclic orientations of $G$.  So if $G$ has a loop then
it has no acyclic orientations and $\cA(G)=\emp$.  A {\it sink} of a
digraph is a vertex $v_0$ such that all arcs incident with $v_0$ are
directed  towards it.  Let $\cA(G,v_0)$ be the set of acyclic 
orientations of $G$ with a unique sink at $v_0$.

A {\it proper coloring} of $G$ with color set $C$ is a map $\ka:V\ra C$
such that $uv\in E$ implies $\ka(u)\neq\ka(v)$.
Now consider the {\it chromatic polynomial} of $G$ which is
$$
\chi_G(n)=\mbox{\# of proper $\ka:V\ra\{1,2\ldots,n\}$}.
$$
It is well known~\cite{roberts} that $\chi_G(n)$ is a polynomial in
$n$ of degree $d=|V|$ so we write
$$
\chi_G(n)=a_0+a_1 n+\cdots+a_d n^d.
$$
If we need to be specific about the graph, we will write $a_i(G)$ for
the coefficient of $n^i$ in $\chi_G(n)$.

Stanley~\cite{sta:aog} was the first to connect acyclic orientations
of graphs and the characteristic polynomial.  In what follows,
absolute value signs around a set denote its cardinality.
\bth[Stanley]						\label{sta}
For any graph $G$
$$
|\cA(G)|=|\chi_G(-1)|.\quad \Qed
$$
\eth
The result of Greene and Zaslavsky that will interest us can
be seen as an analog of Stanley's Theorem for acyclic orientations
with a unique sink~\cite[Theorem 7.3]{gandz}    .
       \begin{theorem}[Greene-Zaslavsky]
						\label{grz}
Let $v_{0}$ be any vertex of $G$.  Then
\beq						\label{grzeq}
|\cA(G,v_0)|=|a_1|. \quad \Qed
\eeq
     \end{theorem}
Originally this theorem was proved using the theory of hyperplane
arrangements.  The purpose of this paper is to give three other proofs
using different techniques.

In the next section we will give a purely inductive proof.
Stanley~\cite{stan} indicated that such a proof exists and we
provide the details.

In the paper just cited, Stanley introduced a symmetric
function analog of the chromatic polynomial and showed that it counts
 the number of acyclic orientations of $G$ with $j$ sinks,
$1\le j \le d$.  Note that this is not quite the same as counting
those with a given sink.  In Section~\ref{csf} we will show how using
noncommutative variables allows us to generalize the Greene-Zaslavsky
Theorem to the level of symmetric functions.

Our final proof is an algorithmic bijection.  To explain it, we need
to recall
Whitney's Broken Circuit Theorem~\cite{whitney}.  A {\it circuit} in a
graph $G$ will 
be the same as a {\it cycle}, i.e., a closed walk with distinct vertices.
If we fix a total order on $E(G)$, a {\em broken circuit} is
a circuit with its largest edge (with respect to the total order)
removed.  Let the {\em broken circuit complex}  $B_{G}$ of $G$ denote
the set of all $S\subseteq E(G)$ which do {\em not} contain a
broken circuit.  The Broken Circuit Theorem asserts:

     \begin{theorem}[Whitney]  \label{nbc}
For any finite graph, $G$, on d vertices we have  
$$
\chr = \sum_{S\in B_{G}}(-1)^{|S|}n^{d-|S|}. 
\quad \Qed
$$
      \end{theorem}

It follows immediately from Theorems~\ref{sta} and~\ref{nbc} that
$|\cA(G)|=|B_G|$.  This result was given a bijective algorithmic
proof by Blass and Sagan~\cite{bs}.  It is also clear from the
previous theorem that 
\beq		\label{a1}
|a_1|=|\{S\in B_G\ :\ |S|=d-1\}|.
\eeq
So to prove the Greene-Zaslavsky Theorem bijectively it suffices to
find a bijection between $\cA(G,v_0)$ and $\{S\in B_G\ :\ |S|=d-1\}$.
This will be done in the last section by modifying the Blass-Sagan
algorithm.

\section{Pure Induction}

We will show that both sides of equation~\ree{grzeq} satisfy the same
recurrence relation and boundary conditions.  We begin with the
well-known Deletion-Contraction Rule for the chromatic
polynomial~\cite{roberts}.
If $e\in E(G)$ we will let $G\setm e$ be $G$ with $e$ deleted.  
We also let $G/e$ be
$G$ with $e$ contracted to a point and any resulting multiple edges
{\it not} identified.  So $|E(G\setm e)|=|E(G/e)|=|E(G)|-1$.  We will
also use this notation for directed graphs.
\bth[Deletion-Contraction Rule]			\label{dc}
For any $e\in E(G)$
$$
\chi_G(n)=\chi_{G\setm e}(n)-\chi_{G/e}(n). \quad\Qed
$$
\eth
{}From this result it is easy to prove inductively that the coefficients
of $\chi_G(n)$ alternate in sign with $a_d=+1$.  Using
Theorem~\ref{dc} again, we see that if $e$ is not a loop then
$$
|a_1(G)|=|a_1(G\setm e)|+|a_1(G/e)|.
$$

We now show bijectively that $|\cA(G,v_0)|$ satisfies the same recursion.
\ble \label{useful1}
Consider any vertex
$v_{0}$, and any edge $e=uv_{0}, u\neq v_{0}$, with the corresponding
arc  $a=\overrightarrow{uv_{0}}$.  The map
$$D \longrightarrow \left \{  \begin{array}{ll}
D\setm a\in {\cal A}(G\setm e,v_{0}) & \mbox{ if $D\setm a
\in{\cal A}(G\setm e,v_{0})$}\\
D/a\in {\cal A}(G/e,v_{0})  & \mbox{ if $D\setm a
\notin{\cal A}(G\setm e,v_{0}),$}
\end{array} \right. $$
 is a bijection between ${\cal A}(G,v_{0})$ and
${\cal A}(G\setm e,v_{0}) \uplus {\cal
A}(G/e,v_{0})$, where  the vertex
of $G/e$ formed by contracting $e$ is labeled
$v_{0}$. \ele
{\bf Proof.}
%
%
We  first prove that this map is
well-defined by showing that in both  cases we
actually obtain an acyclic
orientation with unique sink at $v_{0}$.  This
is clear in the first case by definition.  In
the second, where $D\setm a
\notin{\cal A}(G\setm e,v_{0}),$ it must be true that
$D\setm a$ has sinks both at $u$ and at $v_{0}$
(since deleting a
directed edge of $D$ will not introduce a cycle, nor will it cause us
to lose the sink
at $v_{0}$).   So the orientation $D/a$ will be
in ${\cal A}(G/e,v_{0})$: since $u$ and $v_{0}$
were the only sinks
in $D\setm a$ the contraction  must have a
unique sink at $v_{0}$, and  no new
cycles will be formed.  Hence
this map is well-defined.

To see that this is actually a bijection, we
need only exhibit the inverse.  This  is
obtained by simply orienting
all edges of $G$ as in $D\setm a$ or $D/a$
as appropriate, and then adding in $a$.  It should  be clear
that this  map is also well-defined.
\quad\Qed

For the boundary conditions, we will need the following well-known
result.
\ble					\label{useful2}
If $G$ is connected, then any $D \in {\cal A}(G)$ has at least one
sink.  So if $G$ is arbitrary then for any $D \in {\cal A}(G)$, the
number of sinks is greater than or equal to the number of
components of $G$.  \quad\Qed
\ele

We can now complete the first proof of the Greene-Zaslavsky Theorem by
inducting on the number of non-loops incident with $v_0$.  We have
already verified the recurrence relation, so we need only worry about
the boundary conditions.
If $d=1$, then $${\cal X}_{G}(n) = \left \{
\begin{array}{ll}  n & \mbox{if $G=K_{1},$}\\
       0 & \mbox{if $G$ has loops}.\\
       \end{array} \right.$$
 So in this case, $$ |a_{1}|= \left \{
\begin{array}{ll} 1 & \mbox{ if $G=K_{1},$}\\
     0 & \mbox{ if $G$ has loops}\\
     \end{array} \right \}
     = |{\cal A}(G,\v)|. $$
If $d>1$, then having only loops incident with
$\v$ implies there are least two
components in $G$.  In this
case we can prove inductively from Theorem~\ref{dc} that 
$|a_{1}|=0$ and from Lemma~\ref{useful2} we see
that $|{\cal A}(G,\v)|=0$ as
well.  Thus the boundary conditions match and we are done. \quad\Qed

\section{Chromatic Symmetric Functions} 	\label{csf}

Using his symmetric function generalization, $X_G$, of the chromatic 
polynomial, Stanley~\cite{stan} proved a result related to, but not
quite implying, the one of Greene and Zaslavsky. (See
Theorem~\ref{stanl} at the end of this section.)  In~\cite{chrom} we
introduced an analogue of $X_G$  using  noncommutative variables.  This
allows us to use deletion-contraction techniques on symmetric functions
to prove a generalization of  Greene-Zaslavsky at this level.

We begin with  some  background on symmetric
functions
in noncommuting variables. Much of this follows
from the work of Doubilet \cite{doubilet}
(although he does not
explicitly mention such functions in his paper)
but they differ from those considered by Gelfand, et. al.~\cite{thibon}. 
These noncommutative symmetric functions will
be indexed by {\em
set} partitions (as opposed to integer partitions in the commutative
case).  

We will write $\pi=B_1/B_2/\ldots/B_k$ to denote a partition of
$[d]:=\{1,2,\ldots,d\}$, i.e., $\uplus_{i=1}^k B_i=[d]$.  The $B_i$
are called {\it blocks}.  The set of all partitions of $[d]$ form a
lattice $\Pi_d$ under the partial order  of refinement.  We will
let $\mt$ denote the meet operation (greatest lower bound) in $\Pi_d$.

Now let
$\bx=\{x_{1},x_{2},x_{3}, \ldots \}$ be a set of
{\em noncommuting} variables.   We define the
{\em noncommutative monomial symmetric function},
$m_{\pi}$, by:
\begin{equation}
 m_\pi=m_{\pi}(\bx)= \sum_{i_{1},i_{2},\ldots ,i_{d}}
x_{i_{1}}x_{i_{2}}\cdots x_{i_{d}},
\label{mono}
\end{equation}
where the sum is over all sequences 
$i_{1},i_{2},\ldots ,i_{d}$ of positive integers
$\bbP$
such that $i_{j}=i_{k}$ if and only if $j$ and
$k$ are  in the same block of $\pi$.  For
example,
$$m_{124/3}=x_{1}x_{1}x_{2}x_{1}
+x_{2}x_{2}x_{1}x_{2}  +x_{1}x_{1}x_{3}x_{1}+
x_{3}x_{3}x_{1}x_{3} +\cdots  $$ is
the monomial symmetric function in noncommuting
variables corresponding to the partition $\pi =
124/3$.  The $m_\pi$ are clearly linearly independent over $\bbC$ and
we call the span of $\{m_\pi\ :\ \pi\in\Pi_d, d\ge0\}$ the
{\it algebra of noncommutative symmetric functions}.

The other basis we will be interested in is given by the 
{\em noncommutative elementary symmetric functions} 
\begin{equation}
e_{\pi}=e_\pi(\bx)=\sum_{\sigma:\sigma \wedge \pi
=\mh}m_{\sigma}= \sum_{i_{1},i_{2},\ldots
,i_{d}} x_{i_{1}}x_{i_{2}}\cdots
x_{i_{d}},
\label{elem1}
\end{equation} 
  where the second sum is over all sequences 
$i_{1},i_{2},\ldots ,i_{d} $ of $\bbP$
such that
$i_{j}\neq
i_{k}$ if $j$ and $k$ are both in the same block
of $\pi$.  As an example
\bea
e_{124/3} &=& x_1x_2x_1x_3+x_1x_2x_2x_3+x_1x_2x_3x_3+x_1x_2x_4x_3+\cdots\\
	&=& m_{13/2/4}+m_{1/23/4}+m_{1/2/34}+m_{1/2/3/4}.
\eea

We now introduce a noncommutative version, $Y_G$, of
Stanley's chromatic symmetric function, $X_G$.  The latter is obtained
from the former merely by letting the variables commute.

\begin{definition} For any multigraph $G$
with vertices labeled  $ v_{1}, v_{2}, \ldots,
v_{d}  $  {\em in a fixed order}, define
     \[ Y_{G} = Y_G(\bx)= \su , \]
where the sum is over all proper colorings $\ka:V\ra\bbP$
of $G$.
\end{definition}
As an example, if we let $P_3$ be the path with edge set
$E=\{v_1v_2,v_2v_3\}$ then
\begin{eqnarray*}Y_{P_{3}}&=&
x_{1}x_{2}x_{1}+x_{2}x_{1}x_{2}+x_{1}x_{3}x_{1}
+ \cdots + x_{1}x_{2}x_{3}
+x_{1}x_{3}x_{2} +\cdots +x_{3}x_{2}x_{1} +
\cdots \\
&=&m_{13/2}+m_{1/2/3}.
\end{eqnarray*}
Note that if we let $1^n$ denote the substitution
$x_1=x_2=\cdots=x_n=1$ and $x_i=0$ for $i>n$ then 
$$
X_G(1^n)=Y_G(1^n)=\chi_G(n)
$$
since the only terms
surviving in the sum are those using the first $n$ colors.

We will need two properties of $Y_{G}$;
proofs can be found in~\cite{chrom}.
For the first one, consider $\delta \in {\cal S}_{d}$, the symmetric
group on $[d]$.
We let $\delta$ act on the
vertices of $G$ by
$\delta(v_{i})=v_{\delta(i)}$.  This induces an
action on graphs, denoted $\delta(G)=H$, where
$H$ is just a relabeling of $G$.  We also have an action on
noncommutative symmetric functions given by linearly extending
\[\delta \comp (x_{i_{1}}x_{i_{2}} \cdots 
x_{i_{k}}) \defn 
x_{i_{\delta^{-1}(1)}}x_{i_{\delta^{-1}(2)}}
\cdots 
x_{i_{\delta^{-1}(k)}}.\]
These two actions are compatible.
\bpr
{\bf (Relabeling Proposition)}  For any finite
multigraph $G$, we have
$$\delta \comp Y_{G} = Y_{\delta (G)},$$ where
the vertex order $v_{1},v_{2},\ldots,v_{d}$ is
used in both $Y_{G}$ and
$Y_{\delta(G)}$. \label{relabel}
\epr

In order to allow us to state the Deletion-Contraction Rule for
$Y_G$, we make the following definition.
 \begin{definition} Define an operation
{\em induction},\hspace{.02 in} $\ind$,
on monomials in noncommuting variables
by \[  (x_{i_{1}}x_{i_{2}}  \cdots x_{i_{d-2}}
x_{i_{d-1}})\ind
      ~ = ~ x_{i_{1}}x_{i_{2}}  \cdots
x_{i_{d-2}} x_{i_{d-1}}^{2}\]
and extend linearly.
\end{definition}
{}From equation (\ref{mono}) it
is easy to see that if $\pi\in\Pi_{d-1}$, then 
$m_{\pi}\ind = m_{\pi+(d)}$  where $\pi+(d)\in \Pi_d$ is
the partition obtained from $\pi$ by inserting $d$ into the block with
$d-1$.

\bpr {\bf(Deletion-Contraction Rule)} 
If $e=v_{d-1}v_{d}$ is in $E(G)$ then
$$\non = \nond \: -  \nonc \ind,$$
where the 
contraction of $e=v_{d-1}v_{d}$ is labeled
$v_{d-1}$.
\epr
To illustrate, consider $P_3$ again.  So $e=v_2v_3$ and 
$$Y_{P_{3}}=Y_{P_{2}\uplus\{v_{3}\}} -
Y_{P_{2}}\ind.$$  We then compute
 $$\begin{array}{ll}
Y_{P_{2}\uplus\{v_{3}\}} &=m_{1/2/3}+m_{1/23}
+m_{13/2},\\
Y_{P_{2}}&= m_{1/2},\\
Y_{P_{2}}\ind & =m_{1/2}\ind=m_{1/23}.\\
\end{array}$$
So
\begin{eqnarray*}Y_{P_{3}}&=&m_{1/2/3}+m_{1/23}+
m_{13/2}-m_{1/23}\\
   &=&m_{1/2/3}+m_{13/2},
          \end{eqnarray*}

\bth \label{help2}
Let $\non =\displaystyle{ \sum_{{ \pi \in
\Pi_{d}}}c_{\pi}e_{\pi}}$. Then for any fixed
vertex,
$v_{0}$,
$$
|\cA(G,v_0)|=(d-1)!c_{[d]}.
$$
\eth
 {\bf Proof.}  We  induct on the number of
non-loops in $E$.  If all the
edges of $G$ are
loops, then $$Y_{G}= \left \{ \begin{array}{ll}
e_{1/2/\ldots/d} &  \mbox{ if $G$ has no
edges}\\
                                    0 & \mbox{
if $G$ has loops.}\\
                                    \end{array}
                                    \right.$$
 So $$c_{[d]}=\left \{ \begin{array}{ll} 1 &
\mbox{ if $G=K_{1}$}\\
                                        0 &
\mbox{ if $d>1$ or $G$ has loops}\\
                                       
\end{array}
                                        \right
\} = |{\cal A}(G,\v)|.$$                                   

Now suppose that $G$ has non-loops. Then by the Relabeling Proposition,
we
may choose
$e=v_{d-1}v_{d}$ and 
$\non = \nond - \nonc \ind$.  We are 
only interested in the leading coefficient,
so let
\[ \non=ae_{[d]} + \sum_{{\sigma <
[d]}}a_{\sigma}e_{\sigma},\]
\[\nond=be_{[d]} + \sum_{{\sigma <
[d]}}b_{\sigma}e_{\sigma},\] and
\[\nonc=ce_{[d-1]} + \sum_{{\sigma < [d-1]}}
c_{\sigma}e_{\sigma}\]
where $\leq$ is the partial order on set
partitions.
Using induction and Lemma \ref{useful1}, it suffices
to prove  that
$(d-1)!a=(d-1)!b+(d-2)!c.$

{}From the change of basis formulae  found
in \cite{doubilet}
one gets
\begin{equation}
 e_{\pi} \ind = \sum_{{ \sigma \leq \pi}} {{\mu
(\mh,
\sigma)} \over{ \mu (\mh, \sigma
+(d))}}\sum_{{\tau
\leq
\sigma +(d)}} \mu(\tau, \sigma +(d)) e_{\tau}.
\label{changeup}
\end{equation} 
This permits us to compute the coefficient of
$e_{[d]}$ in $\nonc \ind$. The only term which
contributes comes from $ce_{[d-1]} \ind$, and
 \begin{eqnarray*}
ce_{[d-1]} \ind&=& c\sum_{\sigma \in \Pi_{d-1}} 
{{\mu
(\mh, \sigma)} \over{ \mu (\mh, \sigma
+(d))}}\sum_{{\tau
\leq \sigma +(d)}} \mu(\tau, \sigma +(d))
e_{\tau}\\
 &=& c{{\mu(\mh, [d-1])}\over{\mu (\mh,
[d])}}e_{[d]} +
\sum_{\tau <[d]} d_{\tau}e_{\tau}\\
 &=&{-c \over{d-1}}e_{[d]}+\sum_{\tau <[d]}
d_{\tau}e_{\tau}\\
\end{eqnarray*}
Now $\non = \nond - \nonc \ind$ yields
\begin{eqnarray*} (d-1)!a&=&(d-1)!b+(d-1)!{c
\over{d-1}}\\
 &=&(d-1)!b+(d-2)!c,
\end{eqnarray*}
completing  the proof.    \quad\Qed

To see why this implies Greene-Zaslavsky, 
recall that $Y_G(1^n)=X_{G}(1^{n})={\cal X}_{G}(n)$.  Now, if
$\pi=B_1/B_2/\ldots/B_k$ then under this substitution
$$
e_\pi(1^n)=\prod_{i=1}^k n(n-1)(n-2)\cdots(n-|B_i|+1).
$$
For $k\ge2$, this polynomial is divisible by $n^2$.  So the only 
contribution to the linear term of $\chi_G(n)$ is when $\pi=[d]$ which
yields a coefficient with absolute value $(d-1)!c_{[d]}$.

Immediately from the previous theorem we get
\bco	\label{one}
If $\non=\displaystyle{ \sum_{{ \pi \in
\Pi_{d}}}c_{\pi}e_{\pi}}$, then the number of
acyclic
orientations of
$G$ with one sink is $d!c_{[d]}$.
\eco

Before ending this section, we should state Stanley's
theorem~\cite{stan} relating 
$X_G$ and sinks.  In it, $e_\lambda$ is the commutative elementary
symmetric function corresponding to the integer partition $\lambda$,
and $l(\lambda)$ is the number of parts of $\lambda$.
\bth[Stanley]	\label{stanl}
If $X_{G}=\sum_{\lambda}c_{\lambda}e_{\lambda}$,
then the number of acyclic orientations of $G$
with $j$ sinks is
$\ds{\sum_{l(\lambda)=j}c_{\lambda}.}$
\quad  \Qed
\eth
We can prove an analogue of this this theorem
in the noncommutative setting by using
his technique involving $P$-partitions.
However, this only implies Corollary~\ref{one} (and Theorem~\ref{sta})
but not Theorem~\ref{help2}.

\section{The Modified Blass-Sagan Algorithm}

We will now prove Theorem~\ref{grz} a third time, using~\ree{a1} to
interpret the linear coefficient of $\chi_G(n)$.  This demonstration
will use a variant of an algorithmic bijection of Blass and Sagan
to show that $\cA(G,v_0)$ and $\{S\in B_G\ :\ |S|=d-1\}$ have the same
cardinality.

We first need some notation and definitions. 
For any arc $a=\overrightarrow{wu}$, the
oppositely oriented arc is denoted $a'= \overrightarrow{uw}.$  We
also say that to {\em unorient} an arc, $a$, in
a digraph we will just
add the oppositely oriented arc $a'$.  By the same token,
an edge will also be considered as a pair of oppositely
oriented arcs so that any graph is also a digraph.  Since we
are interested in acyclic digraphs, it is necessary
to adopt the convention
that a digraph is acyclic if it has no cycles
of length $\geq 3$.   With this convention, unorienting an arc will 
not necessarily produce
a cycle.   Also for any acyclic digraph $D$, we
will let $c(D)$ be the
{\em contraction of D},  which is the graph where
all {\em unoriented}
arcs of $D$ have been contracted.  We  note that  $c(D)$ is
still acyclic
and has no unoriented arcs.

\bth
For any fixed vertex $v_{0}\in V(G)$, the number
of
acyclic orientations of $G$ with a unique sink 
at $v_{0}$
is the same as  the number of sets, $S\in B_{G}$
with
$|S|=d-1$.
\eth

{\bf Proof.}
We will construct a bijection using an algorithm  that
sequentially examines each arc of an element of
$\cA(G,v_0)$  and either deletes the arc
or unorients it.  

Fix an orientation of $G$ (not necessarily acyclic) 
which we will refer to as the {\em normal
orientation}, and
also choose a fixed vertex $v_{0}$ of $G$. The
algorithm
will accept any
acyclic orientation $D$ of $G$ which has a
unique sink at
$v_{0}$, and consider
each arc in turn, using the total order on the
edges which defines
the broken
circuits.  At the stage when an arc
$a=\overrightarrow{wu}$ is
being considered,
the algorithm will delete $a$
if either

I) $D\cup a'$ has a cycle, or

II) $c(D)\setm a$ has only one sink, {\em and} $a$ is 
{\em not} normally oriented.

\n Otherwise, the algorithm will unorient $a$.
For an example of how this algorithm works,  see
Figure \ref{alg}.
The steps of the algorithm are
labeled either {\bf I}, {\bf II}, or {\bf u},
indicating if the 
algorithm deleted  the arc for reason I or II,
or unoriented it.

\unitlength=.75mm
\special{em:linewidth 0.4pt}
\linethickness{0.4pt}
\bfi
\begin{picture}(180,160)(0,10)

\put(75,155){\circle*{1.5}}
\put(90,140){\circle*{1.5}}
\put(90,170){\circle*{1.5}}
\put(105,155){\circle*{1.5}}
\put(20,155){\makebox(0,0){Normal Orientation:}}
\put(89,169){\vector(-1,-1){13}}
\put(89,141){\vector(-1,1){13}}
\put(90,169){\vector(0,1){-28}}
\put(104,156){\vector(-1,1){13}}
\put(91,141){\vector(1,1){13}}
\put(80,165){\makebox(0,0){$a_{1}$}}
\put(100,165){\makebox(0,0){$a_{2}$}}
\put(95,155){\makebox(0,0){$a_{4}$}}
\put(100,143){\makebox(0,0){$a_{3}$}}
\put(80,144){\makebox(0,0){$a_{5}$}}
\put(-15,115){\makebox(0,0){$D$}}
\put(0,115){\makebox(0,0){$v_{0}$}}
\put(15,125){\circle*{1.5}}
\put(5,115){\circle*{1.5}}
\put(25,115){\circle*{1.5}}
\put(15,105){\circle*{1.5}}
\put(14,124){\vector(-1,-1){8}}
\put(14,106){\vector(-1,1){8}}
\put(15,124){\vector(0,1){-18}}
\put(16,124){\vector(1,-1){8}}
\put(24,114){\vector(-1,-1){8}}
\put(7,123){\makebox(0,0){$a_{1}$}}
\put(-15,90){\makebox(0,0){$c(D)\setminus a$}}
\put(0,90){\makebox(0,0){$v_{0}$}}
\put(15,100){\circle*{1.5}}
\put(5,90){\circle*{1.5}}
\put(25,90){\circle*{1.5}}
\put(15,80){\circle*{1.5}}
%
\put(14,81){\vector(-1,1){8}}
\put(15,99){\vector(0,1){-18}}
\put(16,99){\vector(1,-1){8}}
\put(24,89){\vector(-1,-1){8}}
\put(30,102.5){\vector(1,0){20}}
\put(40,107){\makebox(0,0){{\bf I}}}
\put(65,125){\circle*{1.5}}
\put(55,115){\circle*{1.5}}
\put(75,115){\circle*{1.5}}
\put(65,105){\circle*{1.5}}
\put(64,106){\vector(-1,1){8}}
\put(65,124){\vector(0,1){-18}}
\put(66,124){\vector(1,-1){8}}
\put(74,114){\vector(-1,-1){8}}
\put(73,122){\makebox(0,0){$a_{2}$}}
\put(65,100){\circle*{1.5}}
\put(55,90){\circle*{1.5}}
\put(75,90){\circle*{1.5}}
\put(65,80){\circle*{1.5}}
\put(64,81){\vector(-1,1){8}}
\put(65,99){\vector(0,1){-18}}
\put(74,89){\vector(-1,-1){8}}
\put(80,102.5){\vector(1,0){20}}
\put(90,107){\makebox(0,0){{\bf II}}}
\put(115,125){\circle*{1.5}}
\put(105,115){\circle*{1.5}}
\put(125,115){\circle*{1.5}}
\put(115,105){\circle*{1.5}}
\put(114,106){\vector(-1,1){8}}
\put(115,124){\vector(0,1){-18}}
\put(124,114){\vector(-1,-1){8}}
\put(123,107){\makebox(0,0){$a_{3}$}}
\put(115,100){\circle*{1.5}}
\put(105,90){\circle*{1.5}}
\put(125,90){\circle*{1.5}}
\put(115,80){\circle*{1.5}}
\put(114,81){\vector(-1,1){8}}
\put(115,99){\vector(0,1){-18}}
\put(30,32.5){\vector(1,0){20}}
\put(40,37){\makebox(0,0){{\bf u}}}
\put(-15,45){\makebox(0,0){$D$}}
\put(50,45){\makebox(0,0){$v_{0}$}}
\put(65,55){\circle*{1.5}}
\put(55,45){\circle*{1.5}}
\put(75,45){\circle*{1.5}}
\put(65,35){\circle*{1.5}}
\put(64,36){\vector(-1,1){8}}
\put(65,54){\vector(0,1){-18}}
\put(74,44){\line(-1,-1){8}}
\put(68,45){\makebox(0,0){$a_{4}$}}
\put(-15,20){\makebox(0,0){$c(D)\setminus a$}}
\put(50,20){\makebox(0,0){$v_{0}$}}
\put(65,30){\circle*{1.5}}
\put(55,20){\circle*{1.5}}
\put(65,10){\circle*{1.5}}
\put(64,11){\vector(-1,1){8}}
\put(80,32.5){\vector(1,0){20}}
\put(90,37){\makebox(0,0){{\bf u}}}
\put(115,55){\circle*{1.5}}
\put(105,45){\circle*{1.5}}
\put(125,45){\circle*{1.5}}
\put(115,35){\circle*{1.5}}
\put(114,36){\vector(-1,1){8}}
\put(115,54){\line(0,1){-18}}
\put(124,44){\line(-1,-1){8}}
\put(107,37){\makebox(0,0){$a_{5}$}}
\put(105,20){\circle*{1.5}}
\put(115,10){\circle*{1.5}}
%
\put(130,32.5){\vector(1,0){20}}
\put(140,37){\makebox(0,0){{\bf u}}}
\put(165,55){\circle*{1.5}}
\put(155,45){\circle*{1.5}}
\put(175,45){\circle*{1.5}}
\put(165,35){\circle*{1.5}}
\put(164,36){\line(-1,1){8}}
\put(165,54){\line(0,1){-18}}
\put(174,44){\line(-1,-1){8}}
\epi
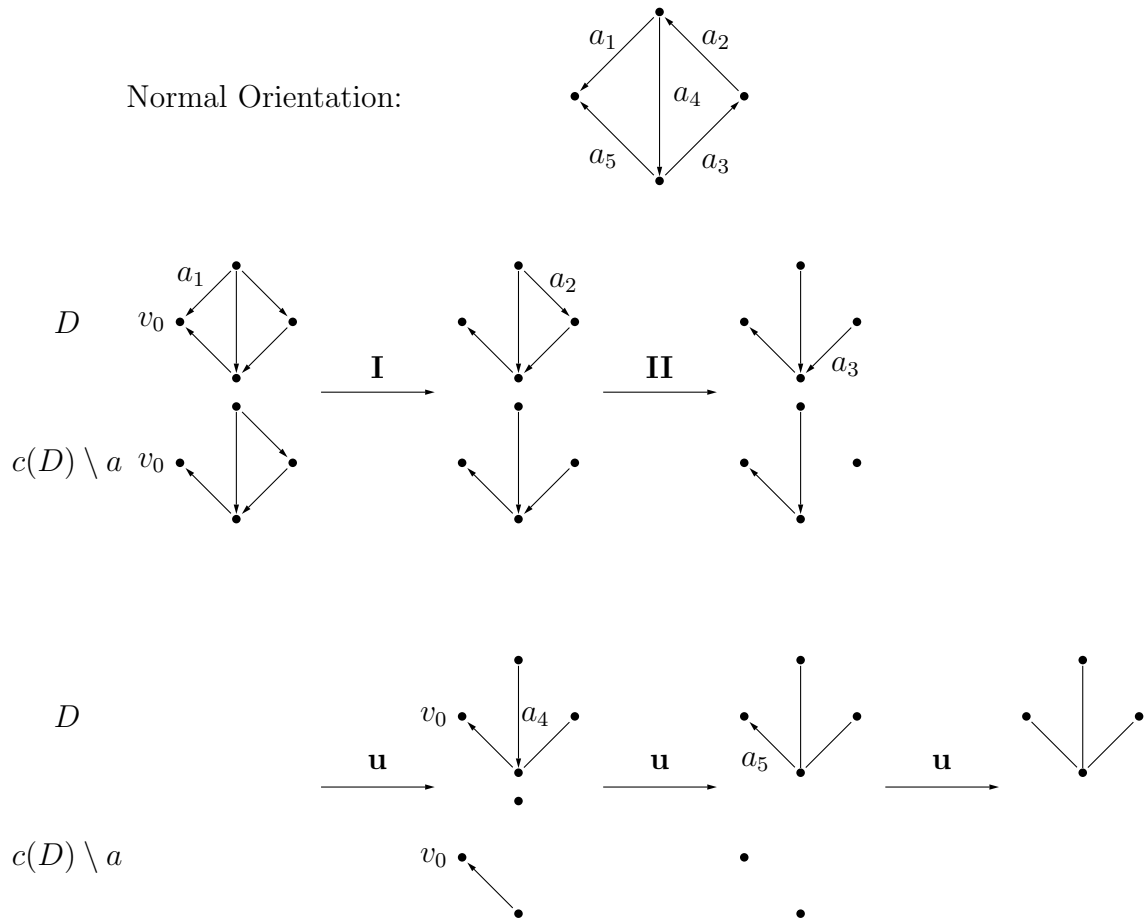
 \capt{An example of the Algorithm} \label{alg}
\efi
To show that this algorithm actually does
produce a
bijection, we shall
first introduce a sequence of sets, 
$\cD_{0},\cD_{1},
\ldots,\cD_{q}$
such that $\cD_{0}$ is the set of all acyclic
orientations of $G$ with a
unique sink at $v_{0}$, and $\cD_{q}$ (where $q =
|E(G)|$)
 is the set of all $S\in B_{G}$ with $|S|=d-1$. 
Equivalently, $\cD_{q}$ is the set of all
spanning trees, $T$, of
$G$ such
that $E(T)$ contains no broken circuits.

We will show  that the $k$th step of the 
algorithm
gives a bijection,
 $f_{k}:\cD_{k-1}\ra\cD_{k}$, where
 $\cD_{k}$  is defined as the
  set of all spanning subdigraphs $D$ of $G$
satisfying the following conditions:

(a) Each of the first $k$ edges of $G$ is either
present in $D$ (as an
unoriented edge) or absent from $D$, but each of
the
remaining $q-k$
edges is present in $D$ in exactly one
orientation.

(b) $D$ is acyclic.

(c)  $D$ has a $x \rightarrow v_{0}$ path for
every $x\in
V(D)$.

(d) The unoriented part of $D$ contains no
broken
circuit.

{}From these conditions, it should be clear that
$\cD_{0}$ is indeed the
set of acyclic orientations of $G$ with  a
unique sink
at $v_{0}$ by Lemma~\ref{useful2}.  It is 
also clear that any element of $\cD_{q}$ will be
an
acyclic, connected
graph, which implies  that the elements of 
$\cD_{q}$
must be trees with exactly
$d-1$ edges.  So provided the algorithm gives
a bijection at each
step, we will have the desired bijection between
acyclic orientations
of $G$ with a unique sink at $v_{0}$, and edge
sets of
size $d-1$ which
contain no broken circuits.

We should also note here that conditions (b)
and (c) together
imply that $c(D)$ must have a unique sink which
occurs
at the vertex
identified with $v_{0}$.  That this is the only
possible
sink of $c(D)$ is
clear from condition (c).   We also know that
$v_{0}$ must
be a sink of
$c(D)$, since if it is not, then there is a
vertex $u$
and arc
$a=\overrightarrow{v_{0}u}$ in $c(D)$.  But from
condition (c) there
would have to be
a $u \rightarrow v_{0}$ path in $D$. This 
contradicts the acyclicity of $D$.

To show that the  algorithm does indeed produce
a
bijection at each
step, we use the following three lemmas.  We
also use
the notational
convention that a digraph in $\cD_{k}$ will be
denoted
by $D_{k}$.

\ble \label{arg}
$f_{k}$ maps $\cD_{k-1}$ into $\cD_{k}$.
\ele

{\bf Proof.} We need only  prove that 
properties (a)-(d)
listed previously are still satisfied after the algorithm is
applied at the $k$th
stage.  We proceed to verify each one in turn.

(a)  Since at the $k$th step the algorithm will
either
delete or
unorient the $k$th arc, this is clear.

(b) Since any arc which would form a cycle if
unoriented will be
deleted by the algorithm, this also is clear.

(c) Since unorienting an arc can never destroy
an $x
\rightarrow v_{0}$
 path, we need only consider the case where the
algorithm  deletes an
arc.  In fact, if the arc 
$a=\overrightarrow{wu}$ in $D_{k-1}$
was deleted, we
 need only show that there is still a $w
\rightarrow
v_{0}$ path.

Now, if the arc $a = \overrightarrow{wu}$ in
$D_{k-1}$ was
deleted for the first
reason, then we must have had another
(different) $w\rightarrow u$ path in $D_{k-1}$. 
 Since
there was a $u \rightarrow v_{0}$ path in
$D_{k-1}$, (in
fact, one which
didn't use the arc $a$)   we can then extend
our other $w\rightarrow u$ path into a
walk containing a $w \rightarrow v_{0}$ path in
$D_{k}$. 

If the arc $a = \overrightarrow{wu}$ in
$D_{k-1}$ was deleted for
the second
reason, again we need only consider the
possibility
that for the vertex
$w$, there is no $w \rightarrow v_{0}$ path in
$D_{k}$. 
But then there is no {\em oriented} arc
$\overrightarrow{wv}$ with $u\neq v$, since
 otherwise all $v\rightarrow
v_{0}$ paths must also use the arc  $a$, as
there are no $w \rightarrow v_{0}$ paths in
$D_{k}$. Thus
$D_{k-1}$ would have a cycle containing $w$. 
Contracting 
all 
unoriented arcs from $w$ and repeating this
argument as
necessary, we see
 that $w$
would then be a sink of $c(D_{k-1})\setm a$, which 
contradicts our reason
for deleting $a$.

(d)  Suppose for the sake of contradiction that
the
unoriented part of
$D_{k}$ contains a broken circuit, $C \setm x$, where
$x$
is the greatest
element of the cycle $C$.  Since the unoriented
part
of $D_{k-1}$ didn't
 contain any broken circuits, and since the only
difference between
$D_{k-1}$ and $D_{k}$ is at the $k$th arc $a$,
we see
that $a$ must be
unoriented in $D_{k}$ and that $a\in C\setm x$.  But
then
$x$ is greater than
$a$, and so $x$ is present in $D_{k}$ in one of
its
orientations.  But
all the other edges in $C$ are also present and
unoriented.  Hence, $C$
forms a cycle in $D_{k}$, contradicting the
previously
verified fact that $D_{k}$ is acyclic.   
\quad\Qed

\ble
$f_{k}$ is one-to-one.
\ele

{\bf Proof.}  Suppose $D_{k-1}$ and $D_{k-1}'$ are
two
distinct elements of
$\cD_{k-1}$ which are both mapped to $D_k$ by the
algorithm.  Since the
algorithm only affects the $k$th arc, we note
that 
$D_{k-1}$ and
 $D_{k-1}'$ (and consequently $c(D_{k-1})$ and
$c(D_{k-1}')$)
must only differ
 in that arc.  Without loss of generality,
we may assume that this arc is $a$
with normal orientation in $D_{k-1}$ and $a'$ with
abnormal orientation in $D_{k-1}'$.

We note that  $D_k$ was not obtained from $D_{k-1}$ and
$D_{k-1}'$ by 
deletion.  For  if $a$ was
deleted
from $D_{k-1}$ for the first reason then $D_{k-1}'$ has a
cycle and vice-versa.  And the
second reason does not apply to $a$ which has normal orientation.

If the $k$th arc was unoriented then, by reason II,
$c(D_{k-1}')\setm a'$ must have an
additional sink.  So if $a'=\overrightarrow{uw}$
then $u$ must be the
extra sink.  But this
means that $u$ is also an additional sink in $c(D_{k-1})$,
contradicting $D_{k-1}\in
\cD_{k-1}$.  \quad\Qed

\ble
$f_{k}$ maps $\cD_{k-1}$ onto $\cD_{k}$.
\ele

{\bf Proof.}  Given $D_{k}\in \cD_{k}$ we must
construct $D_{k-1}
\in \cD_{k-1}$ which maps onto it. Hence for any
digraph,  $D_{k}\in
\cD_{k}$, we must construct a digraph $D_{k-1}$
and
verify that the
algorithm does indeed map $D_{k-1}$ onto
$D_{k}$, and
that $D_{k-1}$
satisfies properties (a)-(d).  For all of the
following cases, it will
 be immediate that the $D_{k-1}$ we construct
will
satisfy properties
(a), (b), and (d), so we will only do the
verification of property (c).
Let $e$ be the $k$th edge of $G$. There are two
cases.

The first case is  when $e$ is not an edge of $D_{k}$.  If 
there
exists a unique
orientation $a$ of $e$ in which $D_{k}$ would remain
acyclic, we give $e$
that orientation in $D_{k-1}$.  If both
orientations
of $e$ would
preserve the acyclicity of $D_{k}$, then we
choose $a$ to be the
abnormal
orientation for $e$ in $D_{k-1}$.  We note that
at
least one of the
orientations of $e$ must preserve acyclicity,
since
otherwise $e$
completes two different cycles in $D_{k-1}$. 
These two cycles
  together would contain a cycle in $D_{k}$,
which is
a contradiction.

That the algorithm maps the digraph $D_{k-1}$
obtained
in the previous
paragraph to $D_{k}$ is obvious when only one
orientation of $a$ produces
an acyclic orientation of $D_{k-1}$.  However,
if both
produce acyclic
orientations, we need to check that
$c(D_{k-1})\setm a$ has
a unique sink at
$v_{0}$. This is true, since it is easy to see
that
$c(D_{k-1})\setm a=c(D_{k-1}\setm a)= c(D_{k})$.  To
verify that
$c(D_{k-1})$
constructed above still satisfies property (c),
we
note that adding
 an arc cannot destroy any existing paths. 
So the first case is done. 

In the second case we have $e$ present in $D_k$
and so  neither orientation can
produce a cycle
in $D_{k-1}$. We note that there must be at
least one
orientation of $e=wu$ such that there remains an
$x
\rightarrow v_{0}$ path for every $x\in
D_{k-1}$.  If all
$x \rightarrow v_{0}$ paths $P$ use the arc
$a=\overrightarrow{wu}$
for some $x$, and if all $y \rightarrow v_{0}$
paths $Q$
use $a'=\overrightarrow{uw}$ for some $y$, then the
$x\rightarrow w$ portion
of $P$ together with the $w\rightarrow v_{0}$
portion of
$Q$ contains an $x \rightarrow v_{0}$ path
avoiding $a$,
which contradicts our assumption about $x$.

 If there is a unique orientation of $e=wu$ so
that
there
 remains an $x \rightarrow v_{0}$ path for every
$x\in
D_{k-1}$  we choose
that one to maintain property (c) for $D_{k-1}$,
say $a=\overrightarrow{wu}$. Using the same
argument we used to
prove the second case of (c) in Lemma \ref{arg},
it
is easy to verify
that the algorithm will
take the $D_{k-1}$ so constructed and map it to
$D_{k}$ by unorienting
$a$ since 
$c(D_{k-1})\setm a$ has an additional sink at $w$ .  

In the subcase where $e$ is present in
$D_{k}$ as
an unoriented edge
and we  would still retain property (c) with
either
orientation of $e$,
  we will consider the digraph $D_{k-1}$
obtained from
$D$ by giving $e$ the
normal orientation, say $a=\overrightarrow{wu}$.
  It is clear
that the algorithm maps $D_{k-1}$ to $D_{k}$,
since
$D_{k-1}\cup a' = D_{k}$ is acyclic and $a$ has the
normal
orientation.   \quad\Qed

\end{document}